\RequirePackage{fix-cm}
\documentclass[12pt]{article}       
\usepackage{lineno,hyperref}
\usepackage{amsmath,mathrsfs,amssymb,enumerate}
\newtheorem{theorem}{Theorem}
\newtheorem{corollary}{Corollary}[theorem]
\newtheorem{conjecture}{Conjecture}[theorem]
\newtheorem{remark}{Remark}[theorem]
\begin{document}

\title{Complete monotonicity of some entropies
}


\author{Ioan Ra\c{s}a\\ Department of Mathematics,
Technical University of Cluj-Napoca,\\
Memorandumului Street 28,\\
400114 Cluj-Napoca,\\
Romania, ioan.rasa$@$math.utcluj.ro
}
\date{}

\maketitle

\begin{abstract}
It is well-known that the Shannon entropies of some parameterized probability distributions are concave functions with respect to the parameter. In this paper we consider a family of such distributions (including the binomial, Poisson, and negative binomial distributions) and investigate the Shannon, R\'{e}nyi, and Tsallis entropies of them with respect to the complete monotonicity.\\
\textbf{keywords}: entropies; concavity; complete monotonicity; inequalities\\
\textbf{subject class}: 94A17; 60E15; 26A51
\end{abstract}

\section{Introduction}
\label{intro}

Let $c\in \mathbb{R}$, $I_c := \left [ 0, -\frac{1}{c}\right ]$ if $c<0$, and $I_c:= [0,+\infty )$ if $c \geq 0$.

For $\alpha \in \mathbb{R}$ and $k \in \mathbb{N}_0$ the binomial coefficients are defined as usual by
\begin{equation*}
{\alpha \choose k}:=\frac{\alpha (\alpha -1)\dots (\alpha-k+1)}{k!}\quad \text{if } k \in \mathbb{N}, \text{ and } {\alpha \choose 0}:=1.
\end{equation*}

Let $n> 0$ be a real number such that $n>c$ if $c\geq 0$, or $n=-cl$ with some $l\in \mathbb{N}$ if $c<0$.

For $k\in \mathbb{N}_0$ and $x\in I_c$ define
\begin{equation*}
p_{n,k}^{[c]}(x):=(-1)^k {-\frac{n}{c} \choose k}(cx)^k (1+cx)^{-\frac{n}{c}-k}, \quad \text{if } c\neq 0,
\end{equation*}
\begin{equation*}
p_{n,k}^{[0]}(x):=\lim _{c\to 0} p_{n,k}^{[c]}(x)= \frac{(nx)^k}{k!}e^{-nx}.
\end{equation*}

Details and historical notes concerning these functions can be found in~\cite{3}, \cite{7}, \cite{21} and the references therein. In particular,
\begin{equation}
\frac{d}{dx}p_{n,k}^{[c]}(x) = n \left ( p_{n+c,k-1}^{[c]}(x) - p_{n+c,k}^{[c]}(x)\right ).\label{eq:1}
\end{equation}

Moreover,
\begin{equation}
\sum _{k=0}^\infty p_{n,k}^{[c]}(x) = 1;\label{eq:2}
\end{equation}

\begin{equation}
\sum _{k=0}^\infty k p ^{[c]}_{n,k}(x)=nx,\label{eq:3}
\end{equation}
so that $\left (p_{n,k}^{[c]}(x)\right )_{k\geq 0}$ is a parameterized probability distribution. Its associated Shannon entropy is
\begin{equation*}
H_{n,c}(x):=-\sum_{k=0}^\infty p_{n,k}^{[c]}(x) \log p_{n,k}^{[c]}(x),
\end{equation*}
while the R\'{e}nyi entropy of order $2$ and the Tsallis entropy of order $2$ are given, respectively, by (see~\cite{18}, \cite{20})
\begin{equation*}
R_{n,c}(x):= -\log S_{n,c}(x); \quad T_{n,c}(x):=1-S_{n,c}(x),
\end{equation*}
where
\begin{equation*}
S_{n,c}(x) := \sum _{k=0}^\infty \left (p_{n,k}^{[c]}(x)\right )^2, \quad x\in I_c.
\end{equation*}

The cases $c=-1$, $c=0$, $c=1$ correspond, respectively, to the binomial, Poisson, and negative binomial distributions. For other details see also~\cite{15}, \cite{16}.

In this paper we investigate the above entropies with respect to the complete monotonicity.

\section{Shannon entropy}

\subsection*{A. Let's start with the case $c<0$.}

$H_{n,-1}$ is a concave function; this is a special case of the results of~\cite{19}; see also~\cite{6}, \cite{8}, \cite{9} and the references therein.

Here we shall determine the signs of all the derivatives of $H_{n,c}$.

\begin{theorem}
Let $c<0$. Then, for all $k\geq 0$,
\begin{equation}
H_{n,c}^{(2k+2)}(x)\leq 0, \quad x \in \left ( 0,-\frac{1}{c} \right ),\label{eq:4}
\end{equation}

\begin{equation} H_{n,c}^{(2k+1)}(x) = \label{eq:5}
 \begin{cases}
    \geq 0          & x \in  ( 0,-\frac{1}{2c}  ],\\
    \leq 0          & x \in [ -\frac{1}{2c}, - \frac{1}{c}  ).\\
  \end{cases}
\end{equation}
\end{theorem}

\textbf{Proof}
We have $n=-cl$ with $l \in \mathbb{N}$. As in~\cite{10}, let us represent $\log{(l!)}$ by integrals:
\begin{equation}
\log{(l!)} = \int _0 ^\infty \left ( l - \frac{1-e^{-ls}}{1-e^{-s}} \right )\frac{e^{-s}}{s} ds = \int _0 ^1 \left ( \frac{1-(1-t)^l}{t} -l \right ) \frac{dt}{\log{(1-t)}}.\label{eq:6}
\end{equation}

Now using~\eqref{eq:2}, \eqref{eq:3} and \eqref{eq:6} we get

\begin{equation*}
H_{n,c}(x) = H_{l,-1}(-cx) = - l \left [(-cx)\log{(-cx)}+(1+cx)\log{(1+cx)}\right ]+
\end{equation*}
\begin{equation*}
\int _0 ^1 \frac{-t}{\log{(1-t)}} \frac{(1+cxt)^l+(1-t-cxt)^l-1-(1-t)^l}{t^2}dt.
\end{equation*}

It is a matter of calculus to prove that

\begin{eqnarray*}
H''_{n,c}(x) &=& cl \left ( \frac{1}{x} - \frac{c}{1+cx}\right )  \\&+& c^2l(l-1)\int _0 ^1 \frac{-t}{\log{(1-t)}} \left [ (1+cxt)^{l-2} + (1-t-cxt)^{l-2}\right ] dt,
\end{eqnarray*}
and for $k\geq 0$

\begin{eqnarray*}
&&H_{n,c}^{(2k+2)}(x)=cl(2k)! \left ( \frac{1}{x^{2k+1}} - \left ( \frac{c}{1+cx}\right )^{2k+1} \right )\\
&+& l(l-1)\dots (l-2k-1)c^{2k+2}\\&& \int _0 ^1 \frac{-t}{\log{(1-t)}} \left [ (1+cxt)^{l-2k-2} + (1-t-cxt)^{l-2k-2}\right ]t^{2k} dt.
\end{eqnarray*}

For $0<t<1$ we have
\begin{equation}
0<\frac{-t}{\log{(1-t)}}<1, \label{eq:new7}
\end{equation}
so that
\begin{equation}
H_{n,c}^{(2k+2)}(x) \leq cl(2k)! \left ( \frac{1}{x^{2k+1}} - \left ( \frac{c}{1+cx}\right )^{2k+1} \right )+\label{eq:7}
\end{equation}
\begin{equation*}
+ l(l-1)\dots (l-2k-1)c^{2k+2} \int _0 ^1  \left [ (1+cxt)^{l-2k-2} + (1-t-cxt)^{l-2k-2}\right ]t^{2k} dt.
\end{equation*}

Repeated integration by parts yields
\begin{equation*}
\int _0 ^1 (1+cxt)^{l-2k-2}t^{2k}dt \leq \frac{(2k)!}{(l-2)(l-3)\dots (l-2k-1)(cx)^{2k}}\int _0 ^1 (1+cxt)^{l-2}dt,
\end{equation*}
and so
\begin{equation}
\int _0 ^1 (1+cxt)^{l-2k-2}t^{2k}dt \leq \frac{(2k)!\left [ (1+cx)^{l-1}-1 \right ]}{(l-1)(l-2)\dots (l-2k-1)(cx)^{2k+1}}.\label{eq:8}
\end{equation}

Replacing $x$ by $-\frac{1}{c}-x$ we obtain
\begin{equation}
\int _0 ^1 (1-t-cxt)^{l-2k-2}t^{2k}dt \leq \frac{(2k)! \left [ 1-(-cx)^{l-1} \right ]}{(l-1)(l-2)\dots (l-2k-1)(1+cx)^{2k+1}}.\label{eq:9}
\end{equation}

From~\eqref{eq:7}, \eqref{eq:8} and \eqref{eq:9} it follows that
\begin{equation*}
H_{n,c}^{(2k+2)}(x)\leq cl(2k)! \left [ \frac{(1+cx)^{l-1}}{x^{2k+1}} - \frac{c^{2k+1}(-cx)^{l-1}}{(1+cx)^{2k+1}}\right ] \leq 0,
\end{equation*}
and this proves~\eqref{eq:4}.

It is easy to verify that $H_{n,c}^{(2k+1)}\left ( -\frac{1}{2c} \right ) = 0$. Since $H_{n,c}^{(2k+2)}\leq 0$, it follows that $H_{n,c}^{(2k+1)}$ is decreasing, and this implies~\eqref{eq:5}.

\subsection*{B. Consider the case $c=0$.}

$H_{n,0}$ is the Shannon entropy of the Poisson distribution. The derivative of this function is completely monotonic: see, e.g., \cite[p. 2305]{2}. For the sake of completeness we insert here a short proof.

\begin{theorem}
$H'_{n,0}$ is completely monotonic, i.e.,
\begin{equation}
(-1)^k H_{n,0}^{(k+1)}(x)\geq 0, \quad k \geq 0,\quad x>0. \label{eq:10}
\end{equation}
\end{theorem}

\textbf{Proof}
Let us remark that $H_{n,0}(y) = H_{1,0}(ny)$; so it suffices to investigate the derivatives of $H_{1,0}(x)$.

According to \cite[(2.5)]{10},
\begin{eqnarray*}
H_{1,0}(x) &=& x-x\log {x} + \int _0 ^\infty \frac{e^{-t}}{t} \left (x - \frac{1-\exp {(x(e^{-t}-1))}}{1-e^{-t}} \right )dt\\
&=& x-x\log{x} - \int _0 ^1 \left ( x - \frac{1-e^{-sx}}{s} \right ) \frac{ds}{\log{(1-s)}}.
\end{eqnarray*}

It follows that
\begin{equation*}
H'_{1,0}(x) = -\log{x} - \int _0^1 \left ( 1-e^{-sx}\right )\frac{ds}{\log{(1-s)}}
\end{equation*}
and for $k\geq 1$,
\begin{equation}
H_{1,0}^{(k+1)}(x) = (-1)^k \left ( \frac{(k-1)!}{x^k} + \int _0^1 s^k e^{-sx} \frac{ds}{\log{(1-s)}}\right ). \label{eq:11}
\end{equation}

By using~\eqref{eq:new7} we get
\begin{equation*}
\int _0 ^1 \frac{s^k e^{-sx}}{\log{(1-s)}}ds \geq - \int _0 ^1 s^{k-1}e^{-sx} ds =
\end{equation*}

\begin{equation*}
=-\int _0 ^x \frac{t^{k-1}}{x^k} e^{-t}dt \geq - \int _0 ^\infty \frac{1}{x^k}t^{k-1}e^{-t}dt = -\frac{(k-1)!}{x^k}.
\end{equation*}

Combined with~\eqref{eq:11}, this proves~\eqref{eq:10} for $k\geq 1$. In particular, we see that $H_{n,0}$ is concave and non-negative on $[0,+\infty )$; it follows that $H'_{n,0}\geq 0$ and so~\eqref{eq:10} is completely proved.

\subsection*{C. Let now $c>0$.}

\begin{theorem}
For $c>0$, $H'_{n,c}$ is completely monotonic.
\end{theorem}

\textbf{Proof}
Since $H_{m,c}(y) = H_{\frac{m}{c},1}(cy)$, it suffices to study the derivatives of $H_{n,1}(x)$.

By using~\eqref{eq:2}, \eqref{eq:3} and
\begin{equation*}
\log{A} = \int _0 ^\infty \frac{e^{-x}-e^{-Ax}}{x}dx, \quad A>0,
\end{equation*}
we get
\begin{equation*}
H_{n,1}(x) = n \left ( (1+x)\log{(1+x)} - x\log{x} \right ) +\int _0 ^\infty \frac{e^{-ns}-e^{-s}}{s(1-e^{-s})}\left ( 1-(1+x-xe^{-s})^{-n} \right )ds
\end{equation*}
\begin{equation*}
 = n \left ( (1+x)\log{(1+x)} - x\log{x} \right ) + \int _0 ^1 \frac{1-(1-t)^{n-1}}{t\log{(1-t)}} \left (1-(1+tx)^{-n} \right ) dt.
\end{equation*}

It follows that, for $j\geq 1$,
\begin{equation*}
\frac{1}{n}H_{n,1}^{(j+1)}(x) = (-1)^{j-1}(j-1)! \left ( (x+1)^{-j} - x^{-j} \right )+
\end{equation*}
\begin{equation*}
+(-1)^{j-1}(n+1)(n+2)\dots (n+j) \int _0 ^1 \frac{-t}{\log{(1-t)}} \left [ 1-(1-t)^{n-1} \right ] (1+xt)^{-n-j-1}t^{j-1}dt.
\end{equation*}

Using again~\eqref{eq:new7}, we get
\begin{eqnarray*}
&&(-1)^{j-1}\frac{1}{n}H_{n,1}^{(j+1)}(x) \leq (j-1)! \left ( (x+1)^{-j}-x^{-j}\right )+\\
&&+(n+1)(n+2)\dots (n+j) \int _0 ^1 \left [1-(1-t)^{n-1}\right ] (1+xt)^{-n-j-1}t^{j-1}dt\\
&&=u(x) + v(x),
\end{eqnarray*}
where
\begin{equation*}
u(x):=\frac{(j-1)!}{(x+1)^j} - (n+1)(n+2)\dots (n+j) \int _0 ^1 t^{j-1}(1-t)^{n-1} (1+xt)^{-n-j-1}dt,
\end{equation*}
\begin{equation*}
v(x):= (n+1)(n+2)\dots (n+j) \int _0 ^1 t^{j-1} (1+xt)^{-n-j-1}dt - \frac{(j-1)!}{x^j}.
\end{equation*}

We shall prove that $u(x)\leq 0$ and $v(x)\leq 0$, $x>0$. Let us remark that
\begin{equation}
\int _0 ^1 t^{j-1}(1-t)^{n-1}(1+xt)^{-n-j-1}dt \geq \int _0 ^1 t^{j-1}(1-t)^n(1+xt)^{-n-j-1}dt, \label{eq:12}
\end{equation}
and integration by parts yields
\begin{equation*}
\int _0 ^1 \frac{t^{j-1} (1-t)^n}{(1+xt)^{n+j+1}}dt = \frac{j-1}{(n+1)(x+1)} \int _0 ^1 \frac{t^{j-2}(1-t)^{n+1}}{(1+xt)^{n+j+1}}dt.
\end{equation*}

Applying repeatedly this formula we obtain
\begin{equation}
\int _0 ^1 \frac{t^{j-1}(1-t)^n}{(1+xt)^{n+j+1}}dt = \frac{(j-1)!}{(n+1)(n+2)\dots (n+j)}\frac{1}{(x+1)^j}.\label{eq:13}
\end{equation}

Now \eqref{eq:12} and \eqref{eq:13} imply $u(x) \leq 0$.

Using again integration by parts we get
\begin{eqnarray*}
\int _0 ^1 t^{j-1}(1+xt)^{-n-j-1}dt \leq \frac{j-1}{(n+j)x} \int _0 ^1 t^{j-2}(1+xt)^{-n-j}dt \\ \leq \dots \leq \frac{(j-1)!}{(n+1)(n+2)\dots (n+j)}\frac{1}{x^j},
\end{eqnarray*}
which shows that $v(x) \leq 0$.

We conclude that
\begin{equation}
(-1)^{j-1}H_{n,1}^{(j+1)}(x)\leq 0, \quad j \geq 1, x>0. \label{eq:14}
\end{equation}

In particular, \eqref{eq:14} shows that $H_{n,1}$ is concave on $[0,+\infty )$; it is also non-negative, which means that $H'_{n,1}\geq 0$. Combined with \eqref{eq:14}, this shows that $H'_{n,1}$ is completely monotonic, and the proof is finished.

\begin{remark}
\eqref{eq:13} can be obtained alternatively by using the change of variables $y=(1-t)/(1+xt)$ and the properties of the Beta function. An alternative proof of the inequality $v(x)\leq 0$ follows from
\begin{equation*}
\int _0^1 t^{j-1} (1+xt)^{-n-j-1} dt \leq \frac{1}{x^{j-1}}\int _0 ^\infty \frac{(xt)^{j-1}}{(1+xt)^{n+j+1}}dt=
\end{equation*}
\begin{equation*}
=\frac{1}{x^j} \int _0 ^\infty \frac{s^{j-1}}{(1+s)^{j+n+1}}ds=\frac{1}{x^j}B(j,n+1)=\frac{1}{x^j}\frac{(j-1)!n!}{(n+j)!}.
\end{equation*}

\end{remark}

\begin{corollary}
The following inequalities are valid for $x>0$ and $c\geq 0$:
\begin{equation}
\log{\frac{x}{cx+1}} \leq \sum _{k=0}^\infty p_{n+c,k}^{[c]}(x)\log{\frac{k+1}{ck+n}}\leq \log{\frac{nx+1}{ncx+n}}.\label{eq:15}
\end{equation}

In particular, for $c=0$ and $n=1$,
\begin{equation*}
\log{x} \leq \sum _{k=0}^\infty e^{-x}\frac{x^k}{k!}\log{(k+1)}\leq \log{(x+1)}.
\end{equation*}
\end{corollary}

\textbf{Proof}
We have seen that $H'_{n,c}(x)\geq 0$. An application of \eqref{eq:1} yields
\begin{equation*}
H'_{n,c}(x) = n \left ( \log{\frac{1+cx}{x}} + \sum _{k=0}^\infty p_{n+c,k}^{[c]}(x)\log{\frac{k+1}{n+ck}} \right ).
\end{equation*}

This proves the first inequality in~\eqref{eq:15}; the second is a consequence of Jensen's inequality applied to the concave function $\log{t}$.

\section{R\'{e}nyi entropy and Tsallis entropy}

The following conjecture was formulated in~\cite{13}:
\begin{conjecture}\label{conj:3.1}
$S_{n,-1}$ is convex on $[0,1]$.
\end{conjecture}

Th. Neuschel~\cite{11} proved that $S_{n,-1}$ is decreasing on $\left [ 0, \frac{1}{2}\right ]$ and increasing on $\left [ \frac{1}{2}, 1\right ]$. The conjecture and the result of Neuschel can be found also in~\cite{5}.

A proof of the conjecture was given by G. Nikolov~\cite{12}, who related it with some new inequalities involving Legendre polynomials. Another proof can be found in~\cite{4}.

Using the important results of Elena Berdysheva~\cite{3}, the following extension was obtained in~\cite{17}:
\begin{theorem}\label{th:3.2}
(\cite[Theorem 9]{17}). For $c<0$, $S_{n,c}$ is convex on $\left [ 0, -\frac{1}{c}\right ]$.
\end{theorem}

A stronger conjecture was formulated in~\cite{14} and~\cite{17}:
\begin{conjecture}\label{conj:3.3}
For $c\in \mathbb{R}$, $S_{n,c}$ is logarithmically convex, i.e., $\log S_{n,c}$ is convex.
\end{conjecture}

It was validated for $c\geq 0$ by U. Abel, W. Gawronski and Th. Neuschel~\cite{1}, who proved a stronger result:
\begin{theorem}\label{th:3.4}
(\cite{1}). For $c\geq 0$, the function $S_{n,c}$ is completely monotonic, i.e.,
\begin{equation*}
(-1)^m \left ( \frac{d}{dx} \right )^m S_{n,c}(x)>0, \quad x\geq 0, m\geq 0.
\end{equation*}

Consequently, for $c\geq 0$, $S_{n,c}$ is logarithmically convex, and hence convex.
\end{theorem}

Summing up, for the R\'{e}nyi entropy $R_{n,c} = -\log S_{n,c}$ and Tsallis entropy $T_{n,c}=1-S_{n,c}$, we can state
\begin{corollary}
\begin{enumerate}[i)]
\item{} Let $c\geq 0$. Then $R_{n,c}$ is increasing and concave, while $T'_{n,c}$ is completely monotonic on $[0,+\infty )$.

\item{}$T_{n,c}$ is concave for all $c\in \mathbb{R}$.
\end{enumerate}
\end{corollary}

\textbf{Proof}
\begin{enumerate}[i)]
\item{}Apply Theorem~\ref{th:3.4}.

\item{}For $c<0$, apply Theorem~\ref{th:3.2}. For $c\geq 0$, Theorem~\ref{th:3.4} shows that $S_{n,c}$ is convex, so that $T_{n,c}$ is concave.
\end{enumerate}

\begin{remark}
As far as we know, Conjecture~\ref{conj:3.3} is still open for $c<0$, so that the concavity of $R_{n,c}$, $c<0$, remains to be investigated.
\end{remark}

\subsection*{Acknowledgement}
The author is grateful to the referee for valuable comments and very constructive suggestions. In particular, the elegant alternative proofs presented in Remark 1 were kindly suggested by the referee.



\end{document}